\documentclass{article}
\usepackage{maa-monthly}

\usepackage[T1]{fontenc}
\usepackage{graphicx}     
\usepackage{hyperref} 
\usepackage{url}
\usepackage{amsfonts} 

\usepackage{multicol}
\usepackage{tikz}

\newtheorem{lemma}{Lemma}
\newtheorem{proposition}{Proposition}

\theoremstyle{definition}
\newtheorem*{definition}{Definition}

\newtheorem*{example}{Example}

\begin{document}

\title{Extremal Values of Pi\footnote{This is an Original Manuscript of an article published by Taylor \& Francis in \textit{The American Mathematical Monthly} on 26 Sept.~2022, available at \href{https://maa.tandfonline.com/doi/pdf/10.1080/00029890.2022.2115779}{https://maa.tandfonline.com/doi/pdf/10.1080/00029890.2022.2115779}.\newline Its 2020 \textit{Mathematics Subject Classifications} are: Primary 52A40; Secondary 52A10, 52A38.}}
\markright{Extremal Values of Pi}
\author{Nikhil Henry Bukowski Sahoo}

\maketitle

\begin{abstract}
We discuss the classical results of Stanis\l aw Go\l\k{a}b, on the values of pi in arbitrary normed planes, including the classification of extremal values. We reprove the result of \cite{pi-values}, which states that any norm with quarter-turn symmetry have pi-value $\geq \pi.$ We also show that a norm is Euclidean iff it has quarter-turn symmetry in some basis and has pi-value $\pi$.
\end{abstract}

\noindent In 1932, the Polish geometer Stanis\l aw Go\l\k{a}b posed and solved an interesting problem: if we generalize the notion of $\text{pi}=\text{circumference}/\text{diameter}$ to unit circles of arbitrary norms on the plane, then the possible pi-values comprise the interval $[3,4]$. The article \cite{pi-values} provides a wonderful exposition and proves some related results, including a new theorem showing that norms on $\mathbb R^2$ with quarter-turn symmetry only attain pi-values in the interval $[\pi,4].$ We will recount these results from scratch, while also providing classifications of the norms that achieve extreme pi-values. For the values of $3$ and $4$, this classification was also given by Go\l\k{a}b in his original paper. For the value of $\pi$, with quarter-turn symmetry, this classification is (to the best of my knowledge) new. We will then relate this final result to a fundamental question of Minkowski geometry: what are necessary and sufficient conditions for a normed space to be Euclidean? 

Recall that a norm $X$ on $\mathbb R^n$ is given by a function $||\cdot||_X:\mathbb R^n\rightarrow [0,\infty)$ such that: \begin{enumerate}
    \item[$\bullet$] $||v||_X=0$ if and only if $v=0$;
    \item[$\bullet$] $||cv||_X=|c|\cdot ||v||_X$ for all $c\in \mathbb R$;
    \item[$\bullet$] and $||u+v||_X\leq ||u||_X+||v||_X.$
\end{enumerate} 
Typical examples include the $\ell^p$ norms on $\mathbb R^n$ for all $p\geq 1$:
\begin{equation}\label{def-lp}
    ||v||_p=\sqrt[p\hspace*{-0.4mm}]{{|v_1|}^p+\dots+{|v_n|}^p}.
\end{equation} Taking the limit of (\ref{def-lp}) as $p\rightarrow \infty$ gives the $\ell^\infty$ norm: $||v||_\infty=\max\big\{|v_1|,\dots,|v_n|\big\}.$
The most common norm is $\ell^2$, for which (\ref{def-lp}) is essentially just the ``distance formula.'' For any norm $X$ on $\mathbb R^n$, we may define the unit ball and its boundary, the unit sphere: $$B_X=\{v\in \mathbb R^n:||v||_X\leq 1\}\qquad\text{and}\qquad \partial B_X=\{v\in \mathbb R^n:||v||_X=1\}.$$ Then $B_X$ is compact and convex, with $-B_X=B_X$ and $0\in B_X^\circ$ (the interior of $B_X$). Conversely, for any $B\subset \mathbb R^n$ with these properties and any $v\in \mathbb R^n$, we may define $$||v||_X=\frac{1}{\sup\{a\in \mathbb R:av\in B\}}$$ This is the unique norm $X$ with $B_X=B.$ This bijection between norms and certain convex sets gives the study of normed spaces a geometric flavor: the shape of the unit ball regulates properties of the norm. In particular, we  will see that the aforementioned extreme pi-values occur iff $B_X$ is an ellipse, parallelogram, or affine regular hexagon.

In what follows, we will mostly be focused on two dimensions. As such, we write $\mathcal M$ for the set of norms on $\mathbb R^2.$ It is convenient to identify $\mathbb R^2=\mathbb C$ and use notation befitting complex numbers. In particular, a counter-clockwise quarter-turn is $z\mapsto iz,$ and the $\ell^2$ norm can be written as the absolute value $|z|=||z||_2.$ We speak of angles in terms of the function $\arg:\mathbb C\setminus 0\rightarrow S^1,$ where $S^1$ denotes the circle (the reader may think of $S^1$ as the quotient $\mathbb R/2\pi \mathbb Z$ or as the $\ell^2$ unit circle $\{z\in \mathbb C:|z|=1\}.$) We will always speak of angle-measure in terms of radians.

\section{Arc-length in terms of a norm}~

\medskip\noindent Recall that a metric on a set $A$ is a symmetric function $d:A\times A\rightarrow [0,\infty)$ such that:
\begin{enumerate}
    \item[$\bullet$] $d(u,v)=0$ if and only if $u=v$;
    \item[$\bullet$] and $d(u,w)\leq d(u,v)+d(v,w)$.
\end{enumerate} A norm $X$ on $\mathbb R^n$ induces a metric $d_X(u,v)=||u-v||_X$ on $\mathbb R^n$, with the additional properties $d_X(u+w,v+w)=d_X(u,v)$ and $d(au,av)=a\cdot d(u,v)$ for any $a>0$ and $u,v,w\in \mathbb R^n$ (conversely, any metric on $\mathbb R^n$ with these properties defines a norm). 

A metric $d$ on a set $A$ defines some notion of ``shortest distance'' between two points, but it can also be useful to consider the length travelled along more circuitous paths. Given $u,v\in A,$ a path from $u$ to $v$ is any function $\varphi:[a,b]\rightarrow A$ such that $\varphi(a)=u$ and $\varphi(b)=v.$ (If $u=v$, then we call $\varphi$ a loop with base-point $u.$) The length of $\varphi$ is $$\text{len}_d\varphi=\sup\left\{\sum_{i=1}^nd\big(\varphi(t_{i-1}),\varphi(t_i)\big):a=t_0\leq \dots\leq t_n=b\text{ and }n\in \mathbb N\right\}.$$ We call the sequence $\varphi(t_0),\dots,\varphi(t_n)$ a partition of the path $\varphi$. The resulting length always takes a well-defined value in $[0,\infty].$ First, we note some elementary properties. As the results are intuitively plausible, the proofs are left to the interested reader. 
\begin{enumerate}
    \item[(a)]For any paths $\varphi_1:[a,b]\rightarrow A$ and $\varphi_2:[b,c]\rightarrow  A$ such that $\varphi_1(b)=\varphi_2(b)$,\footnote{If two paths share an endpoint, we can always shift their domains, so that the domains line up in this way. By (a) and (c), the way in which this is done will not affect the length of the concatenation.} we define their concatenation $\varphi_1\bullet \varphi_2:[a,c]\rightarrow  A$ as
    $$(\varphi_1\bullet \varphi_2)(t)=\left\{\begin{array}{cc}
         \varphi_1(t),& a\leq t\leq b \\
         \varphi_2(t),& b\leq t\leq c
    \end{array}\right.$$
    In the arithmetic of $[0,\infty]$, we then have $\text{len}_d(\varphi_1\bullet \varphi_2)=\text{len}_d(\varphi_1)+\text{len}_d(\varphi_2).$ 
    
    \item[(b)] Given a loop $\varphi:[a,c]\rightarrow A$ and any $b\in [a,c]$, we may define
    $$\varphi'(t)=\left\{\begin{array}{cc}
         \varphi(t), & b\leq t\leq c \hspace*{13mm} \\
         \varphi(t-c+a), & c\leq t\leq b+c-a 
    \end{array}\right.$$
    This loop $\varphi':[b,b+c-a]\rightarrow A$ is essentially $\varphi$ with a shifted base-point. Using (a), we can see that $\text{len}_d\varphi'=\text{len}_d\varphi$, so the basepoint doesn't matter when measuring length. As such, we will henceforth view loops as functions $S^1\rightarrow A,$ which possess a well-defined notion of length (measured from any basepoint).
    
    \item[(c)] If $\varphi:[c,d]\rightarrow A$ is a path and $f:[a,b]\rightarrow [c,d]$ is monotonic and surjective, then $\text{len}_d(\varphi\circ f)=\text{len}_d\varphi.$ Thus if $\varphi$ is an injective, continuous path (or loop), then $\text{len}_d\varphi$ depends only on the image of $\varphi$. 
    
    \item[(d)] For any $u,v\in \mathbb R^n$, we use the notation $[u,v]=\{(1-t)u+tv:0\leq t\leq 1\}$ and $(u,v)=\{(1-t)u+tv:0<t< 1\}$ for the closed and open segments from $u$ to $v$. For any sequence $x_0,x_1,\dots,x_n\in \mathbb R^n$, we have a polygonal path $$[x_0,x_1,\dots,x_n]=[x_0,x_1]\bullet \cdots\bullet [x_{n-1},x_n].$$
    For any norm $X$ on $\mathbb R^n$, the length in terms of $d_X$ of such a path is 
    $$\text{len}_X[x_0,x_1,\dots,x_n]=d_X(x_0,x_1)+\dots+d_X(x_{n-1},x_n).$$ By additivity of lengths in (a), it suffices to prove that $\text{len}_X[u,v]=||u-v||_X.$
    
    \item[(e)] Suppose that $X$ and $Y$ are norms on $\mathbb R^n$ with $B_X\subset B_Y$. Then $||v ||_Y\leq ||v||_X$ for any $v\in \mathbb R^n$, so we also get $\text{len}_Y\varphi\leq \text{len}_X\varphi$ for any path $\varphi:[a,b]\rightarrow \mathbb R^n.$
\end{enumerate}
Henceforth, we will only ever consider lengths with respect to a norm.

\subsection{Convex paths} In general, curves need not have finite length with respect to a norm, even if they are continuous and injective.\footnote{A typical example is the Koch snowflake. In general, a continuous, injective curve or loop has finite length if and only if its Hausdorff dimension is 1 (or 0, in the case when the ``curve'' is just a single point) \cite{fractals}.} However, we are interested in curves that form a portion of the boundary of a convex set in $\mathbb R^2$. Below, we will prove a useful comparison lemma; in  particular, this will imply that all such curves have finite length.

If $B\subset \mathbb R^n$ is a compact, convex set with $B^\circ\neq \emptyset,$ we will call $B$ a convex body. Then a set $B\subset \mathbb R^n$ is the unit ball of a norm on $\mathbb R^n$ if and only if $B$ is a convex body and symmetric about the origin. If $B\subset \mathbb R^2$ is a convex body, then $\partial B$ is a continuous, injective loop: choose any $x\in B^\circ$ and define $\varphi:\partial B\rightarrow S^1$ by $\varphi(y)=\arg(y-x)$ (since $x\in B^\circ$, we have $y-x\neq 0$). The resulting length depends only on the set $\partial B$, so any $x\in B^\circ$ works equally well. We now compare the lengths of $\partial B$ for various $B$.

\begin{lemma}\label{convex-loops}
If $B_1,B_2\subset \mathbb R^2$ are any two convex bodies with $B_1\subset B_2$, then we have $\textup{len}_X(\partial B_1)\leq \textup{len}_X(\partial B_2)$ for any norm $X\in \mathcal M.$
\end{lemma}

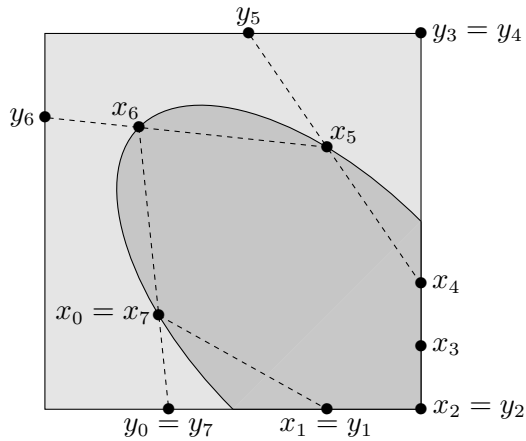
\begin{figure}
    \centering
    \vspace*{5mm}
    \begin{tikzpicture}[scale=2.5]
    \draw[fill=lightgray,fill opacity=0.4] (1,1)--(1,-1)--(-1,-1)--(-1,1)--cycle;
    \draw[rotate=-45,fill=lightgray,fill opacity=0.8] ({cos(45)},{sin(45)}) arc(90:270:{2/sqrt(2)} and {1/sqrt(2)});
    \draw[fill=lightgray,fill opacity=0.75] (1,0) -- (1,-1) -- (0,-1);

    \draw (1,-1) node{$\bullet$};
    
    \draw[dash pattern = on 2 off 2] (1,-{1/3}) node{$\bullet$} -- ({1+4/(1-12*sqrt(1/5))},1) node{$\bullet$};
    
    \draw[dash pattern = on 2 off 2] (-0.5,0.5) node{$\bullet$} -- ({1-3*sqrt(1/5)},-1) node{$\bullet$};
    
    \draw[dash pattern = on 2 off 2] (0.5,{sqrt(4/5)-0.5}) node{$\bullet$} -- (-1,{1-sqrt(1/5)}) node{$\bullet$};
    
    \draw[dash pattern=on 2 off 2] ({-sqrt(4/5)+0.5},-0.5) node{$\bullet$}  -- (0.5,-1) node{$\bullet$};

    \draw (1,1) node{$\bullet$};
    
    \draw (1,-{2/3}) node{$\bullet$};
    
    \draw (-0.69,-0.5) node{$x_0=x_7$};
    
    \draw (0.495,-1.1) node{$x_1=y_1$};
    
    \draw ({1-3*sqrt(1/5)},-1.1) node{$y_0=y_7$};
    
    \draw (-1.12,{1-sqrt(1/5)}) node{$y_6$};
    
    \draw ({1+4/(1-12*sqrt(1/5))},1.09) node{$y_5$};
    
    \draw (1.31,-1) node{$x_2=y_2$};
    
    \draw (1.3,1) node{$y_3=y_4$};

    \draw (1.135,{-2/3}) node{$x_3$};
    \draw (1.135,{-1/3}) node{$x_4$};
    
    \draw (0.59,0.46) node{$x_5$};
    
    \draw (-0.57,0.59) node{$x_6$};

    \end{tikzpicture}
    
    \caption{Comparing the lengths of curves bounding convex regions}
    \label{fig:compare-loops}
\end{figure}

\begin{proof}
This proof follows \cite{pi-values} extremely closely. Let $x_0,\dots,x_n$ be a partition of $\partial B_1.$ Since $\partial B_1$ is a loop, we have $x_0=x_n.$ For each $i=1,\dots,n$, we define
$$y_i=x_{i-1}+\sup\{t\in [0,\infty):x_{i-1}+t(x_i-x_{i-1})\in B_2\}(x_{i}-x_{i-1}).$$ This is the furthest point along the ray $R=\{x_{i-1}+t(x_i-x_{i-1})\in \mathbb R^2:t\in [0,\infty)\}$ that is also contained in $B_2.$ In particular, we see that $y_i\in \partial B_2$ (often, but not always, $y_i$ is the unique point in $R\cap \partial B_2$). We also set $y_0=y_n$, so the sequence $y_0,\dots,y_n$ is a partition of $\partial B_2.$ This is illustrated in Figure \ref{fig:compare-loops}. Because $x_i\in [x_{i-1},y_i]$, we have
$$d_X(x_{i-1},x_i)+d_X(x_i,y_i)=d_X(x_{i-1},y_i)\leq d_X(x_{i-1},y_{i-1})+d_X(y_{i-1},y_i)$$ for all $i=1,\dots,n$, by the triangle inequality. Summing these inequalities, we have
$$\sum_{i=1}^nd_X(x_{i-1},x_i)+\sum_{i=1}^nd_X(x_i,y_i)\leq \sum_{i=1}^nd_X(x_{i-1},y_{i-1})+\sum_{i=1}^nd_X(y_{i-1},y_i).$$
Since $x_0=x_n$ and $y_0=y_n$, the middle two sums are equal and therefore $$\sum_{i=1}^nd_X(x_{i-1},x_n)\leq \sum_{i=1}^nd_X(y_{i-1},y_i)\leq \text{len}_X(\partial B_2).$$  The partition $x_0,\dots,x_n$ was arbitrary, so this gives $\text{len}_X(\partial B_1)\leq \text{len}_X(\partial B_2).$
\end{proof}

If $B\subset \mathbb R^2$ is a convex body, then we have $B\subset [-a,a]^2$ for large enough $a>0$. Polygonal curves have finite length, so Lemma \ref{convex-loops} gives $\text{len}_X(\partial B)<\infty$. Therefore, the boundary of any convex body must have finite length. More generally:

\begin{definition}
We will say that path $\varphi$ from $u$ to $v$ is \textit{convex} if $\varphi=[u,v]$ or if $\varphi \bullet [v,u]$ is the boundary of a convex body. We write $\varphi_1\prec \varphi_2$ if $\varphi_1$ and $\varphi_2$ are both convex paths from $u$ to $v$, such that $\text{conv}(\varphi_1)\subset \text{conv}(\varphi_2)$ ($\text{conv}$ denotes the convex hull).
\end{definition}

Given a convex path $\varphi$, we have $\text{len}_X\varphi\leq \text{len}_X\big(\varphi\bullet [v,u]\big)<\infty$ for any $X\in \mathcal M.$ Thus all convex paths have finite length. We will need some more results concerning convex paths and their lengths (a different proof of (a) occurs in \S\!\S4.3-4.4 of \cite{thompson}).

\begin{lemma}\label{convex-paths}
\begin{enumerate}
    \item[(a)] If $\varphi_1\prec \varphi_2$, then $\textup{len}_X\varphi_1\leq \textup{len}_X\varphi_2$ for any norm $X\in \mathcal M.$
    \item[(b)] If $\varphi$ is a convex path between two distinct points $u,v\in \mathbb R^2$ and $\varphi\cap (u,v)\neq \emptyset,$ then $\varphi=[u,v].$
    \item[(c)] If $B\subset \mathbb R^2$ is a convex body, then any path along its boundary (i.e.~a non-empty, closed, connected subset of $\partial B$) is a convex path.
    \item[(d)] Let $\varphi$ be a convex path and suppose that the points $p,q,r\in \varphi$ occur in this written order along the path $\varphi$. If $q\in (p,r)$, then $[p,r]\subset \varphi.$
\end{enumerate}

\end{lemma}

\begin{proof} (a) We leave the case when $\varphi_1=[u,v]$ to the reader. If $\varphi_2$ is a line segment, then $\varphi_1\prec \varphi_2$ implies that $\varphi_1$ is as well. Hence, for $i=1$ or 2, we see that $\text{conv}(\varphi_i)$ is a convex set bounded by the loop $\partial\,\text{conv}(\varphi_i)=\varphi_i\bullet [v,u].$ Thus, Lemma \ref{convex-loops} gives
\begin{align*}
    \text{len}_X\varphi_2+d_X(v,u)&=\text{len}_X\big(\varphi_2\bullet [v,u]\big)\\
    &\geq \text{len}_X\big(\varphi_1\bullet [v,u]\big)=\text{len}_X\varphi_1+d_X(v,u),
\end{align*}
since $\text{conv}(\varphi_1)\subset \text{conv}(\varphi_2).$ The desired inequality follows immediately.

\smallskip (b) If $\varphi \bullet [v,u]$  is the boundary of a convex body, then $\varphi$ and $[v,u]$ only intersect at $u$ and $v,$ so $\varphi\cap (u,v)=\emptyset.$ This contradicts our initial assumption, so $\varphi=[u,v].$

\smallskip (c) For brevity, we only sketch this proof. Let $\varphi$ be a path from $u$ to $v$ along $\partial B.$ Then either $(u,v)\subset B^\circ$ or $(u,v)\subset \partial B$.\footnote{Let $A\subset \mathbb R^2$ be any convex set. If $x\in A$ and $y\in A^\circ$, then $(x,y)\subset A^\circ$ (the proof is left to the reader). This can then be used to show that, for any $x,z\in A,$ we have either $(x,z)\subset \partial A$ or $(x,z)\subset A^\circ.$} But if $(u,v)\subset \partial B$, then either $\varphi=[u,v]$ or $\varphi\bullet [v,u]=\partial B$, so $\varphi$ is convex. Thus we suppose that $(u,v)\subset B^\circ$. Let $H$ denote the (closed) half-plane containing $\varphi$, which is cut out by the line through $u$ and $v.$ Because $H$ is convex and closed, we can see that $B\cap H$ is a convex body, such that $\partial (B\cap H)=\varphi\bullet [v,u].$ This shows that $\varphi$ is a convex path in any of the above cases.

\smallskip (d) Let $\psi$ be the portion of $\varphi$ between $p$ and $r$. By assumption, we then have $q\in \psi$, and by (c), we see that $\psi$ is a convex path. Since $q\in \psi\cap (p,r)$ and $\psi$ is a convex path from $p$ to $r$, we see that $[p,r]=\psi\subset \varphi$ by (b). \end{proof}

\section{What values does pi take?}~

\medskip\noindent With the above notion of length, we can now define the promised generalization of pi. 

\begin{definition}
For any norm $X\in \mathcal M$, define $\varpi(X)=\text{len}_X(\partial B_X)/2.$ This symbol $\varpi$ was historically used as a cursive $\pi$. The characters mean two different things to us: we write $\varpi$ for ``pi in terms of a norm'' and $\pi=3.14159\dots$ for the classic constant. 
\end{definition}

The results of the previous section show us that $\varpi(X)\in (0,\infty)$ for any norm $X$. Since $-B_X=B_X$, we see that the loop $\partial B_X$ has half-turn symmetry, so we can also calculate $\varpi(X)$ as the length of the intersection of $\partial B_X$ with the upper half-plane. 

In this section, we will eventually show that the image is $\varpi(\mathcal M)=[3,4]$. Already, we can treat an important class of polygonal examples to show that $[3,4]\subset \varpi(\mathcal M)$.

\begin{figure}
    \centering
        \vspace*{5mm}
    \begin{tikzpicture}[scale=1.6]
        \draw[<->,thick] (-1.3,0) -- (1.3,0);
        \draw[<->,thick] (0,-1.2) -- (0,1.4);
        \draw (1,0) -- (0.5,1) -- (-1,1) -- (-1,0) -- (-0.5,-1) -- (1,-1) -- cycle;
        \draw[dash pattern = on 2 off 2] (0.5,0) -- (0.5,1);
        \draw[dash pattern = on 2 off 2] (0,0) -- (-0.5,1);
        
        \draw (0.5,-0.15) node{$t$};
        \draw (-1.05,0.05) -- (-1.15,0.05);
        \draw (-1.05,0.975) -- (-1.15,0.975);
        \draw (-1.1,0.05) -- (-1.1,0.975);
        \draw (-0.975,1.05) -- (-0.975,1.15);
         \draw (0.475,1.05) --(0.475,1.15);
          \draw (-0.975,1.1) -- (0.475,1.1);
          \begin{scope}[shift={(0.15,0.075)}]
              \draw (0.95,0) --  (0.475,0.95);
          \draw (0.9,-0.025) -- (1,0.025);
          \draw (0.425,0.925) -- (0.525,0.975);
          \draw (-1.35,0.43) node{1};
          \draw (0.81,0.55) node{1};
          \draw (-0.5,1.15) node{$1+t$};
          \end{scope}
          
    \end{tikzpicture}
    \caption{Assuming all values $\varpi(X_t)\in [3,4]$}
    \label{fig:values-3-4}
\end{figure}
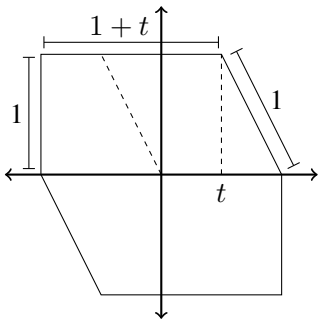

\begin{example}

 Fix $0\leq t\leq 1$ and let $X_t\in \mathcal M$ be the norm whose unit disk is $$\text{conv} \{e_1,te_1+e_2,e_2-e_1,-e_1,-te_1-e_2,e_1-e_2\}.$$
This is shown in Figure \ref{fig:values-3-4} (note that $B_{X_t}$ is a hexagon except in the case of $t=1$, when it degenerates into a square). Since $e_1,e_2,e_2-(1-t)e_1\in \partial B_{X_t},$ we see that these vectors have norm 1 under $X_t.$ Hence, the upper half of the unit circle has length\begin{align*}
    \varpi(X_t)&=\text{len}_{X_t}[e_1,te_1+e_2,e_2-e_1,-e_1]\\&=d_{X_t}(e_1,te_1+e_2)+d_{X_t}(te_1+e_2,e_2-e_1)+d_{X_t}(e_2-e_1,-e_1)\\&=||e_2-(1-t)e_1||_{X_t}+(1+t)||e_1||_{X_t}+||e_2||_{X_t}=3+t.
\end{align*}
This shows that the image of $\varpi:\mathcal  M\rightarrow(0,\infty)$ contains the interval $[3,4],$ as desired.

\end{example}

 Notice that $\varpi(X)$ depends on $X$ in two ways: the unit circle $\partial B_X$ is defined by $X$, but we also use $X$ to measure length. These two dependencies are precisely balanced, in a way that we will now make precise.

\begin{definition}
Let $T:\mathbb R^n\rightarrow \mathbb R^n$ be a linear isomorphism. Given any norm $X$ on $\mathbb R^n$, the push-forward norm $TX$ is defined by $||v||_{TX}=||T^{-1}v||_X$. If norms $X_1$ and $X_2$ satisfy $TX_1=X_2$ for some isomorphism $T$, we will call them linearly equivalent. Similarly, we say that two sets $A_1,A_2\subset\mathbb R^n$ are linearly equivalent if $T(A_1)=A_2$ for some linear isomorphism $T.$
\end{definition}

First, note that $B_{TX}=T(B_X)$ for any norm $X$ on $\mathbb R^n$ and linear isomorphism $T$. Thus we can see that norms $X_1$ and $X_2$ are linearly equivalent if and only if their unit balls $B_{X_1}$ and $B_{X_2}$ are linearly equivalent. This fact will be of frequent use below.  

Fix an isomorphism $T:\mathbb R^n\rightarrow \mathbb R^n$ and a norm $X$ on $\mathbb R^n$. For $u,v\in \mathbb R^n$, we have
$$d_X(u,v)=||u-v||_X=||Tu-Tv||_{TX}=d_{TX}(Tu,Tv).$$
Now consider an arbitrary curve $\varphi:[a,b]\rightarrow \mathbb R^n$.
If $x_0,\dots,x_n$ is a partition of $\varphi$, then $T(x_0),\dots,T(x_n)$ is a partition of $T\circ \varphi$ and we have  $$\sum_{i=1}^nd_X(x_{i-1},x_i)=\sum_{i=1}^nd_{TX}\big(T(x_{i-1},T(x_i)\big)\leq \text{len}_{TX}(T\circ \varphi).$$ This proves that $\text{len}_X\varphi\leq \text{len}_{TX}(T\circ \varphi).$ Replacing $T$ by $T^{-1}$ gives
$$\text{len}_{TX}(T\circ \varphi)\leq \text{len}_{T^{-1}TX}(T^{-1}\circ T\circ \varphi)=\text{len}_X\varphi.$$ Thus $\text{len}_X\varphi=\text{len}_{TX}(T\circ \varphi).$ Note that $\partial B_{TX}=T(\partial B_X)$, so when $n=2$, we have $$\varpi(X)=\text{len}_X(\partial B_X)/2=\text{len}_{TX}(\partial B_{TX})/2=\varpi(TX).$$

\begin{lemma}\label{linear-equiv}
The function $\varpi:\mathcal M\rightarrow (0,\infty)$ is constant on linear equivalence classes.
\end{lemma}

There are three linear equivalence classes of particular importance:
\begin{enumerate}
    \item[(a)] The constant $\pi$ is defined as $\pi=\varpi(\ell^2)$. Because $B_{\ell^2}$ is the classical unit circle, we see that $X\in \mathcal M$ is linearly equivalent to $\ell^2$ if and only if $B_X$ is an ellipse.   Therefore, all ellipses (centered at the origin) yield $\varpi=\pi.$
    
    \item[(b)] Notice that the unit disk that defining $X_1$ in the above example is $$B_{\ell^\infty}=\{(x,y)\in \mathbb R^2:-1\leq x,y\leq 1\}.$$ Hence $\ell^\infty=X_1$ and so $\varpi(\ell^\infty)=\varpi(X_1)=4.$ Since this unit circle is a square, a norm $X\in \mathcal M$ is linearly equivalent to $\ell^\infty$ if and only if $B_X$ is a parallelogram. Therefore, all parallelograms (centered at the origin) yield $\varpi=4.$ In particular, since $B_{\ell^1}$ is the square with vertices $\{e_1,-e_1,e_2,-e_2\},$ we have $\varpi(\ell^1)=4$.
    
    \item[(c)] We will say that $B\subset \mathbb R^2$ is a \textit{linearly regular hexagon} if $B$ is linearly equivalent to a regular hexagon centered at the origin. This is equivalent to the condition: $$B=\text{conv}\{u,v,v-u,-u,-v,u-v\}$$ for some linearly independent $u,v\in \mathbb R^2$ (this is an actual regular hexagon when $|u|=|v|$ and the angle measure between $u$ and $v$  is $\pi/3$). In the above example, we can see that $B_{X_0}$ is linearly regular (take $u=e_1$ and $v=e_2$). We calculated that $\varpi(X_0)=3$, so we see that all linearly regular hexagons yield $\varpi=3.$
\end{enumerate}

In what follows, we will show that $\varpi(\mathcal M)=[3,4]$ and that (b) and (c) characterize the extremal cases. More specifically,
$\varpi(X)=4$ if and only $B_X$ is a parallelogram (centered at the origin), and $\varpi(X)=3$ if and only if $B_X$ is a linearly regular hexagon.

\subsection{Circumscribed parallelograms} We will now address the upper bound $\varpi\leq 4.$ First, we prove that any norm on $\mathbb R^2$ can be put into a particular ``normalized'' form. 

\begin{lemma}\label{parallelogram}
For any $X\in \mathcal M$, there exists some isomorphism $T:\mathbb R^2\rightarrow \mathbb R^2$ such that $||e_1||_{TX}=||e_2||_{TX}=1$ and $||(x,y)||_{TX}\geq \max \big(|x|,|y|\big).$
\end{lemma}

\begin{proof} This proof is different from \cite{pi-values} and more closely follows Theorem 3.2.1 in \cite{thompson}. We will find vectors $u,v\in B_X$ such that $B_X\subset P$, where
$$P=\{su+tv:|s|\leq 1\text{ and }|t|\leq 1\}.$$ Because $u,v\in \partial P$, we must also have $u,v\in \partial B_X$. Therefore $||u||_X=||v||_X=1.$ If $a=\max\big(|s|,|t|\big)$ for some $s,t\in \mathbb R$, then $aB_X\subset aP$ and thus $||su+tv||_X\geq a.$ Then, if we set $T(u)=e_1$ and $T(v)=e_2$, we will get the desired isomorphism $T$. 

We will view $\mathbb R^2$ as the $xy$-plane in $\mathbb R^3$, so that we may consider cross products. Because $B_X$ is compact, we can find two vectors $u,v\in B_X$ that maximize $||u\times v||_2.$ If $B_X\not\subset P$, then there are some $s,t\in \mathbb R$ with $su+tv\in B_X$, but $|s|>1$ or $|t|>1.$ If  $|s|>1$, then we have the following contradiction:
$$||(su+tv)\times v||_2=||su\times v||_2=|s|\cdot ||u\times v||_2>||u\times v||_2.$$ If $|t|>1$, we get a similar contradiction. Thus $B_X\subset P$, completing the proof.
\end{proof}

Using this result, the proof of $\varpi\leq 4$ is almost immediate. We will also classify all $X\in \mathcal M$ satifying $\varpi(X)=4$, following the reproduction of Sch\"affer's proof in \cite{thompson}.

\begin{proposition}\label{upper-bound} For any norm $X\in \mathcal M$, we have $\varpi(X)\leq 4$. Moreover, we have $\varpi(X)=4$ if and only if $B_X$ is a parallelogram (centered at the origin).
\end{proposition}

\begin{figure}
    \centering
        \vspace*{5mm}
    \begin{tikzpicture}[scale=2]
    \draw[thick,<->] (-1.5,0) -- (1.5,0);
    \draw[thick,->] (0,0) -- (0,1.5);
    \draw[fill=lightgray,fill opacity=0.5] (1,0) arc (0:90:1);
    \draw[fill=lightgray,fill opacity=0.5] (0,1) -- ({-1/3},0.95) -- (-0.7,0.7) -- (-0.95,{1/3}) -- (-1,0);
    \draw[draw=none,fill=lightgray,fill opacity=0.5] (1,0) -- (0,1) -- (-1,0);
    \draw (-1,0) -- (-1,1) -- ({sqrt(2)-1},1) -- (1,{sqrt(2)-1}) -- (1,0);
    
    \begin{scope}[shift={({0.05*sqrt(2)},{0.05*sqrt(2)})}]
    \draw ({sqrt(2)-1+0.025},1-0.025) -- (1-0.025,{sqrt(2)-1+0.025});
    \draw ({sqrt(2)-1+0.025-0.025*sqrt(2)},{1-0.025-0.025*sqrt(2)})--
    ({sqrt(2)-1+0.025+0.025*sqrt(2)},{1-0.025+0.025*sqrt(2)});
    \draw ({1-0.025-0.025*sqrt(2)},{sqrt(2)-1+0.025-0.025*sqrt(2)})--
    ({1-0.025+0.025*sqrt(2)},{sqrt(2)-1+0.025+0.025*sqrt(2)});
    
    \end{scope}
    
    \draw (-1.05,0.05) -- (-1.15,0.05);
        \draw (-1.05,0.975) -- (-1.15,0.975);
        \draw (-1.1,0.05) -- (-1.1,0.975);
        
        \draw (1.05,0.05) -- (1.15,0.05);
        \draw (1.05,{sqrt(2)-1.025})-- (1.15,{sqrt(2)-1.025});
        \draw (1.1,{sqrt(2)-1.025}) -- (1.1,0.05);
        
        \draw (1.19,0.23) node{$\xi$};

        \draw (-1.1,0.05) -- (-1.1,0.975);
        \draw (-0.975,1.05) -- (-0.975,1.15);
         \draw ({sqrt(2)-1.025},1.05) --({sqrt(2)-1.025},1.15);
          \draw (-0.975,1.1) -- ({sqrt(2)-1.025},1.1);
         
          \draw (-0.35,1.225) node{$1+\xi$};
        
          \draw (-1.2,0.48) node{1};
        
        \draw (2,0.875) node{$(1-\xi)||e_2-e_1||_X=2(1-\xi)$};
    \end{tikzpicture}\hspace*{-3.5cm}
    \caption{Classifying the case when $\varpi=4$}
    \label{fig:parallelogram}
\end{figure}
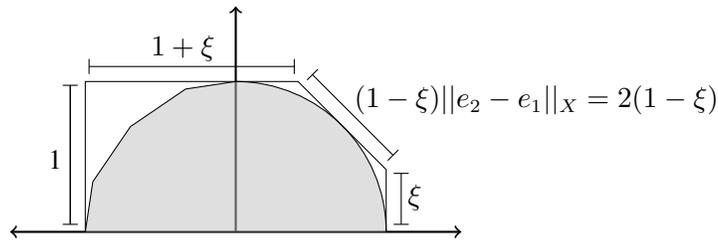

\begin{proof}
Using Lemma \ref{parallelogram}, we may assume that $||(x,y)||_X\geq \max\big(|x|,|y|\big)$ and that $||e_1||_X=||e_2||_X=1$ (since $\varpi$ and the property of being a parallelogram are both preserved under linear equivalence). Then $B_X\subset B_{\ell^\infty}$ and therefore Lemma \ref{convex-loops} gives $$2\varpi(X)=\text{len}_X(\partial B_X)\leq \text{len}_X(\partial B_{\ell^\infty})=4\big(||e_1||_X+||e_2||_X)=8.$$ Now suppose that $\varpi(X)=4.$ If $e_1+e_2,e_1-e_2,-e_1-e_2,e_2-e_1\in B_X,$ then $$B_{\ell^\infty}=\text{conv}\{e_1+e_2,e_1-e_2,-e_1-e_2,e_2-e_1\}\subset B_X\subset B_{\ell^\infty}.$$ Then $X=\ell^\infty$ and $B_X$ is thus a square. Hence, we may assume that one of these points is not in $B_X$; after a rotation, we may assume that $e_1+e_2\notin B_X.$ We define $$\xi=\max\{x+y-1:(x,y)\in B_X\}.$$
Note that $\xi$ exists by the compactness of $B_X$ and $0\leq \xi<1$ because $e_1+e_2\notin B_X.$ Let $\varphi$ be the portion of $\partial B_X$ in the upper half-plane, a path from $e_1$ to $-e_1.$ Then $$\varphi\prec [e_1,e_1+\xi e_2,\xi e_1+e_2,e_2-e_1,-e_1],$$ as illustrated in Figure \ref{fig:parallelogram}. Therefore, Lemma \ref{convex-paths}(a) gives
\begin{align*}
    4=\text{len}_X\varphi&\leq \text{len}_X[e_1,e_1+\xi e_2,\xi e_1+e_2,e_2-e_1,-e_1]\\
    &=||\xi e_2||_X+||(1-\xi)(e_2-e_1)||_X+||(1+\xi)e_1||_X+||e_2||_X\\
    &=2+2\xi+(1-\xi)||e_2-e_1||_X\leq 2(1+\xi)+2(1-\xi)=4,
\end{align*} since $||e_2-e_1||_X\leq ||e_2||_X+||e_1||_X=2.$ We then must have equality throughout; in particular, the last inequality becomes $||e_2-e_1||_X=2$ (since $1-\xi>0$). Hence, $\tfrac{1}{2}(e_2-e_1)\in \partial B_X$ and thus $[e_2,-e_1]\subset \partial B_X$ by Lemma \ref{convex-paths}(d).  Since $e_2-e_1\notin B_X,$ we may repeat this argument with $e_1$ negated, to show that $[e_2,e_1]\subset \partial B_X$ as well. Therefore $\varphi=[e_1,e_2,-e_1]$, which implies that $X=\ell^1$ and $B_X$ is thus a square.
\end{proof}

\subsection{Inscribed hexagons} We now prove that $\varpi\geq 3.$ In classifying the case when $\varpi=3$, we must use the concept of an extreme point, of which we assume no prior knowledge.

\begin{definition}
Let $B\subset \mathbb R^n$ be a convex set. If $F\subset B$ is non-empty, convex and $$x,y\in B\text{ and }(x,y)\cap F\neq \emptyset\quad \Longrightarrow \quad x,y\in F,$$ we call $F$ a face of $B$. We say that $p\in B$ is an extreme point if $\{p\}$ is a face, i.e.
$$x,y\in B\text{ and }p\in (x,y)\quad \Longrightarrow \quad x=y=p.$$
\end{definition}

We will prove that compact, convex sets always contain extreme points. In general, this is proven by induction on dimension (or Zorn's lemma for infinite dimensions), but we only need the case of $\mathbb R^2$. We leave the following facts for the reader to verify:
\begin{enumerate}
    \item[(a)] If $F$ is a face of $B$ and $E\subset F$, then $E$ is a face of $F$ $\Longleftrightarrow $ $E$ is a face of $B$.
    \item[(b)] If $B\subset \mathbb R^n$ is convex and $F\subset B$ is a proper face (i.e.~$F\neq B$), then $F\subset \partial B$.
    \item[(c)] Suppose that $B\subset \mathbb R^n$ is convex and $p\in B$ is the unique point of $B$ maximizing some linear functional $\mathbb R^n\rightarrow \mathbb R.$ Then $p$ is an extreme point of $B$.  
    \item[(d)] For any $p,q\in \mathbb R^n$, the extreme points of the line segment $[p,q]$ are $p$ and $q$.
\end{enumerate}
Our proof that extreme points exist is somewhat odd, but chosen to fit what follows.

\begin{lemma}\label{extreme-pts}\begin{enumerate}
    \item[(a)] Suppose that $A\subset \mathbb R^2$ is convex and $u,v \in A$ are distinct. Let $\ell$ denote the line going through $u$ and $v$. If $[u,v] \subset \partial A$, then $A \cap \ell$ is a face of $A$.
    \item[(b)] If $B\subset \mathbb R^2$ is compact, convex and non-empty, then $B$ has an extreme point.
\end{enumerate}
\end{lemma}

\begin{proof} (a) As an intersection of two convex sets, the set $F=A\cap \ell$ is also convex. Since $u,v\in F,$ we have $F\neq \emptyset$. Suppose that $x,y\in A$ are such that $(x,y)\cap F\neq \emptyset.$ Then the open segment $(x,y)$ must intersect the line $\ell$, so either $x,y\in \ell$ or $x$ and $y$ are (strictly) on opposite sides of $\ell.$ Suppose the latter is true. Let $Q$ be the (possibly non-convex) quadrilateral with boundary $[x,u,y,v,x].$ The convexity of $A$ gives $$Q=\text{conv}\{x,u,v\}\cup \text{conv}\{y,u,v\}\subset A,$$ where we have cut $Q$ along the diagonal $[u,v]$. This diagonal is in the  interior of $Q,$ i.e.~$(u,v)\subset Q^\circ\subset A^\circ$. But this contradicts $[u,v]\subset \partial A,$ so we must have $x,y\in \ell$. Since $x,y\in A$ as well, this shows that $x,y\in F$. Therefore, $F$ is a face of $A$.
 
\smallskip (b) Let $\Pi:\mathbb R^2\rightarrow \mathbb R$ be the projection onto the first coordinate. Define the quantity $m=\max\Pi(B)$ (this exists because $B\neq\emptyset$ is compact) and the line $\ell=\Pi^{-1}(m)$. Note that $B\cap \ell$ is a closed, convex, non-empty subset of $\ell,$ so it is a line segment, i.e.~$B\cap \ell=[u,v]$ with $u,v\in \ell.$ We have $[u,v]=B\cap \ell\subset \partial B$, since otherwise, there would exist $x\in B$ with $\Pi(x)>m.$ Thus if $u\neq v,$ then $[u,v]$ is a face of $B$, by (a). If this is the case, then $u$ is an extreme point of $[u,v]$ and thus of $B$. If $u=v,$ then $\Pi$ is maximized precisely on the set $[u,v]=\{u\},$ so $u$ is an extreme point. \end{proof}

We now prove the bound $\varpi\geq 3$ as in \cite{pi-values}, \cite{golab} or \cite{thompson}, by mimicking the classical straightedge-compass construction of an equilateral triangle. We will also classify all $X\in \mathcal M$ satisfying $\varpi(X)=3$, following the reproduction of Sch\"affer's proof in \cite{thompson}.

\begin{proposition}\label{lower-bound} For any norm $X\in \mathcal M$, we have $\varpi(X)\geq 3$. Moreover, we have $\varpi(X)=3$ if and only if $B_X$ is a linearly regular hexagon.
\end{proposition}

\begin{proof}
By Lemma \ref{extreme-pts}(b), we may choose any extreme point $u\in B_X$ (then $u\in \partial B_X$, so we have $||u||_X=1$). Note that $||0||_X=0$, $||2u||_X=2$ and $0,2u\in \partial B_X+u$. Since $\partial B_X$ is connected, there exists some vector $v\in \partial B_X+u$ with $||v||_X=1.$ Then $v-u\in \partial B_X$, so $||v-u||_X=1.$ This construction is illustrated in Figure \ref{fig:hexagons}(a). For any $0\leq \epsilon<1$, let $w_\epsilon=v+\epsilon u$ and define the hexagon $$H_\epsilon=\text{conv}\{u,w_\epsilon,v-u,-u, -v,u-v\}.$$ 
The vertices are listed in cyclic order (see Figure \ref{fig:hexagons}(a)), which we use to calculate \begin{align*}
\text{len}_X(\partial H_\epsilon) &= ||w_\epsilon-u||_X+||v-u-w_\epsilon||_X\\&\hspace*{2.31cm}+||{-v}||_X+||u-v||_X+||u||_X+||v||_X    \\
&= ||w_\epsilon-u||_X+(1+\epsilon)||u||_X+4=||w_\epsilon-u||_X+5+\epsilon.
\end{align*} 
If $w_\epsilon\in B_X$, then we have $H_\epsilon\subset B_X$ by the convexity of $B_X$, so Lemma \ref{convex-loops} gives
\begin{equation}\label{hexagon-inequality}
    \text{len}_X(\partial B_X)\geq \text{len}_X(\partial H_\epsilon)=||w_\epsilon-u||+5+\epsilon.
\end{equation}
Since $w_0=v\in B_X$ and $||w_0-u||_X=||v-u||_X=1$, this gives $\text{len}_X(\partial B_X)\geq 6.$ This proves that $\varpi(X)\geq 3$, so we now assume that $\varpi(X)=3,$ i.e.~$\text{len}_X(\partial B_X)= 6.$ Notice that $H_0=\text{conv}\{u,v,v-u,-u,-v,u-v\}$ is a linearly regular hexagon. ($H_0$ is also $X\!$-equilateral, meaning that any adjacent vertices are $X\!$-distance 1 apart.) Under the assumption that $\varpi(X)=3$, we will prove that $B_X=H_0.$

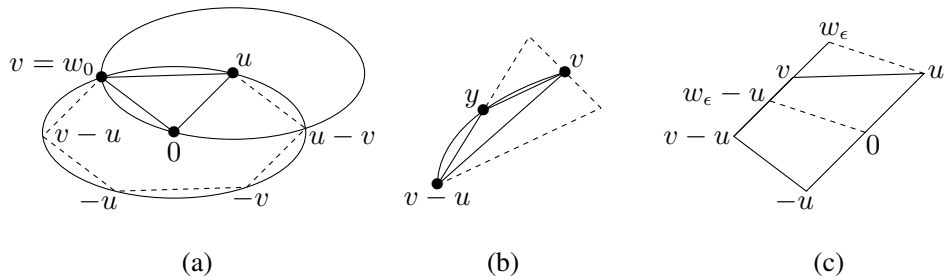
\begin{figure}
    \centering
        \vspace*{5mm}
    \begin{tikzpicture}[scale=1.75]
    
    \coordinate (u) at ({sqrt(1/5)},{sqrt(1/5)});
    \coordinate (nu) at ({-sqrt(1/5)},{-sqrt(1/5)});
    \coordinate (v) at ({1/10*(sqrt(5)-2*sqrt(15))},{1/20*(2*sqrt(5)+sqrt(15))});
    \coordinate (nv) at ({-1/10*(sqrt(5)-2*sqrt(15))},{-1/20*(2*sqrt(5)+sqrt(15))});
    \coordinate (w) at ({sqrt(1/5)-1/10*(sqrt(5)-2*sqrt(15))},{sqrt(1/5)-1/20*(2*sqrt(5)+sqrt(15))});
    \coordinate (nw) at ({-sqrt(1/5)+1/10*(sqrt(5)-2*sqrt(15))},{-sqrt(1/5)+1/20*(2*sqrt(5)+sqrt(15))});

    \begin{scope}[shift={(0,0)}]
    \draw (0.1755,-1) node{(a)};
    
    \draw (0,0) circle (1 and 1/2);
    \draw (u) circle (1 and 1/2);
    \draw[dash pattern=on 2 off 2] (u) -- (w) -- (nv) -- (nu) -- (nw) -- (v);
    
    \draw (u) node{$\bullet$}--(v) node{$\bullet$}--(0,0) node{$\bullet$}--cycle;
    
    \draw (u)++(0.085,0.105) node{$u$};
    \draw (v)++(-0.375,0.08) node{$v=w_0$};
    \draw (nw)++(0.345,0.015) node{$v-u$};
    \draw (w)++(0.28,-0.06) node{$u-v$};
    \draw (0,-0.145) node{$0$};
    \draw (nu)++(-0.13,-0.08) node{$-u$};
    \draw (nv)++(0.03,-0.09) node{$-v$};
    
    \end{scope}

    \draw (2.5,-1) node{(b)};
    
    \begin{scope}[shift={(2,-0.4)}]
    
    \draw (0,0) arc(180:121:2 and 1);
    \draw (0,0) node{$\bullet$} -- (2-1.651,0.5644) node{$\bullet$} -- ({1-0.0301},{0.8572}) node{$\bullet$} -- cycle;
    
    \draw[dash pattern = on 2 off 2](0,0) -- ({2*(1-0.0301+1.651-2)},{2*0.8572-2*0.5644})--(4-2*1.651,2*0.5644)--(2-1.651,0.5644);
    
    \draw ({1-0.0301},{0.8572})++(0.09,0.09) node{$v$};
    
    \draw (2-1.651,0.5644)++(-0.085,0.09) node{$y$};
    
    \draw (0,-0.1) node{$v-u$};
    
    \end{scope}

    \begin{scope}[shift={(5.25,0)}]
    
    \coordinate (u) at ({sqrt(1/5)},{sqrt(1/5)});
    \coordinate (nu) at ({-sqrt(1/5)},{-sqrt(1/5)});
    \coordinate (v) at ({1/10*(sqrt(5)-2*sqrt(15))},{1/20*(2*sqrt(5)+sqrt(15))});
    \coordinate (nv) at ({-1/10*(sqrt(5)-2*sqrt(15))},{-1/20*(2*sqrt(5)+sqrt(15))});
    \coordinate (w) at ({sqrt(1/5)-1/10*(sqrt(5)-2*sqrt(15))},{sqrt(1/5)-1/20*(2*sqrt(5)+sqrt(15))});
    \coordinate (nw) at ({-sqrt(1/5)+1/10*(sqrt(5)-2*sqrt(15))},{-sqrt(1/5)+1/20*(2*sqrt(5)+sqrt(15))});
    
    \draw (-0.2755,-1) node{(c)};
    
    \draw (u) -- (nu) -- (nw) -- (v) -- cycle;
    \draw[dash pattern=on 2 off 2] ({-sqrt(1/5)+1/20*(sqrt(5)-2*sqrt(15))},{sqrt(1/5)-1/40*(2*sqrt(5)+sqrt(15))}) -- (0,0);
    
    \draw[dash pattern=on 2 off 2,shift={(u)}] ({-sqrt(1/5)+1/20*(sqrt(5)-2*sqrt(15))},{sqrt(1/5)-1/40*(2*sqrt(5)+sqrt(15))}) -- (0,0);
    
    \draw[shift={(u)}] (nw) -- ({-sqrt(1/5)+1/20*(sqrt(5)-2*sqrt(15))},{sqrt(1/5)-1/40*(2*sqrt(5)+sqrt(15))});
    
    \draw (u)++(0.09,0.01) node{$u$};
    \draw (nu)++(-0.09,-0.08) node{$-u$};
    \draw (0.06,-0.095) node{$0$};
    \draw (nw)++(-0.28,0) node{$v-u$};
    \draw (v)++(-0.07,0.04) node{$v$};
    
    \draw[shift={(u)}] ({-sqrt(1/5)+1/20*(sqrt(5)-2*sqrt(15))},{sqrt(1/5)-1/40*(2*sqrt(5)+sqrt(15))})++(0.04,0.085) node{$w_\epsilon$};
    
    \draw ({-sqrt(1/5)+1/20*(sqrt(5)-2*sqrt(15))},{sqrt(1/5)-1/40*(2*sqrt(5)+sqrt(15))})++(-0.35,0.05) node{$w_\epsilon-u$};
    
    \end{scope}

    \end{tikzpicture}
    \caption{Three stages in the proof of Proposition \ref{lower-bound}}
    \label{fig:hexagons}
\end{figure}

Let $\varphi$ be the (shorter) path along $\partial B_X$ between any two adjacent vertices of $H_0$. Then we have $\text{len}_X\varphi\geq 1$, since the endpoints of $\varphi$ are $X\!$-distance 1 apart. However, $\text{len}_X(\partial B_X)=6$ is the sum of the lengths of these six paths, so we must have equality $\text{len}_X\varphi=1$ for each path. Now let $\varphi$ be the path along $\partial B_X$ going from $v$ to $v-u$. Since $\text{len}_X\varphi=1$, there is some $y\in \varphi$ with $d_X(v,y)=1/2.$ We also get
$$1=d_X(v,v-u)\leq d_X(v,y)+d_X(y,v-u)\leq \text{len}_X\varphi=1,$$ since $v,y,v-u$ is a partition of $\varphi$. We have equality throughout, so $d_X(v,y)=1/2$ implies that $d_X(y,v-u)=1/2.$ Thus $2(v-y),2(y-v+u)\in \partial B_X$. Notice that
$$u=\tfrac{1}{2}\cdot 2(v-y)+\tfrac{1}{2}\cdot 2(y-v+u)\in \big(2(v-y),2(y-v+u)\big).\footnote{For any vectors $p,q\in \mathbb R^2$, we have $p+q=u$ if and only if $u$ is the midpoint of $[2p,2q]$.}$$ This process is illustrated in Figure \ref{fig:hexagons}(b). Because $u$ is an extreme point, we must have $u=2(v-y)=2(y-v+u)$ and thus $y=v-u/2$ (the midpoint of $[v,v-u]$). Hence $y\in \varphi\cap (v,v-u)$ and therefore $\varphi=[v,v-u]$ by Lemmas \ref{convex-paths}(b) and (c).

If we can prove that $v$ is an extreme point, then we can iterate this whole process, i.e.~show that $[v-u,-u]\subset \partial B_X$, then show that $v-u$ is an extreme point, et cetera. In total, this shows that $\partial B_X=\partial H_0$ and thus $B_X=H_0$ is a linearly regular hexagon.

Consider the line $\ell=\{v+tu:t\in \mathbb R\}$ and note that $v,v-u\in \ell.$ Since we have shown that $[v,v-u]=\varphi\subset \partial B_X$, Lemma  \ref{extreme-pts}(a) states that $B_X\cap \ell$ is a face of $B_X$ and a closed segment (see the proof of said lemma). Thus $v$ is an extreme point of $B_X$ if and only if $v$ is an extreme point of $B_X\cap \ell,$ i.e.~an endpoint of this line segment. But if $v$ is not an endpoint, then there exists $0<\epsilon<1$ with $w_\epsilon=v+\epsilon u\in B_X\cap \ell$. Then we have $w_\epsilon-u=\epsilon v+(1-\epsilon)(v-u)\in [v,v-u]\subset \partial B_X$, so (\ref{hexagon-inequality}) becomes $$6=\text{len}_X(\partial B_X)\geq ||w_\epsilon-u||_X+5+\epsilon=6+\epsilon>6.$$ This is a contradiction, so $v$ must be an endpoint of $B_X\cap\ell$ (see Figure \ref{fig:hexagons}(c)).\end{proof}

\section{Which norms are Euclidean?}~

\medskip\noindent Recall that an inner product is a symmetric, bilinear function $\langle\, \cdot\,,\cdot \,\rangle:\mathbb R^n\times \mathbb R^n\rightarrow \mathbb R$, such that the function $v\mapsto \langle v,v\rangle$ is a norm on $\mathbb R^n.$ If a norm $X$ arises from an inner product in this way, then it is said to be \textit{Euclidean}. The renowned ``parallelogram law'' states that a norm $X$ is Euclidean if and only if, for all $u,v\in \mathbb R^2$, we have
\begin{equation}\label{parallelogram-law}
    2||u||_X^2+2||v||_X^2=||u+v||_X^2+||u-v||_X^2.
\end{equation}
Since (\ref{parallelogram-law}) only involves two vectors at a time, a norm on $\mathbb R^n$ is Euclidean if and only if its restriction to any two-dimensional subspace is Euclidean. This gives a special role to geometric conditions for norms on $\mathbb R^2$ that precisely classify the Euclidean norms.
We will write $\mathcal E\subset \mathcal M$ to denote the set of Euclidean norms on $\mathbb R^2.$

The standard inner product on $\mathbb R^n$ is simply given by $\langle u,v\rangle =u_1v_1+\dots+u_nv_n,$ and it induces the Euclidean $\ell^2$ norm. Any inner product admits an orthonormal basis (using the Gram-Schmidt process), which uniquely characterizes this inner product. Mapping this basis to $e_1,\dots,e_n$, we can see that any Euclidean norm on $\mathbb R^n$ is linearly equivalent to $\ell_2.$ But if $\langle\, \cdot\,,\cdot \,\rangle$ is an inner product and $T:\mathbb R^n\rightarrow \mathbb R^n$ is an isomorphism, then the map $(u,v)\mapsto \langle T^{-1}u,T^{-1}v\rangle$ is also an inner product. Thus Euclidean norms are closed under linear equivalence, so a norm $X$ on $\mathbb R^n$ is Euclidean if and only if $X$ is linearly equivalent to the $\ell^2$ norm (i.e.~$B_X$ is an ellipsoid). 

Returning to two dimensions, we can see that $X\in \mathcal E\Longrightarrow\varpi(X)=\varpi(\ell^2)=\pi$, by Lemma \ref{linear-equiv}. However, the converse does not hold, since we saw an example where $\varpi(X)=\pi$ and $B_X$ is a hexagon. As a first step towards our novel classification of $\mathcal E$, we classify circles via the property ``any tangent is perpendicular to the unique radius that it touches.'' To state this result precisely, we define $\ell_v=v+\text{span}(iv)$. If $v\neq 0$, then $\ell_v$ is the unique line through $v$ that is perpendicular to the line segment $[0,v]$. 

\begin{lemma}\label{classify-circles}
    Suppose that the norm $X\in \mathcal M$ satisfies $B_X^\circ \cap \ell_v=\emptyset$ for all $v\in \partial B_X.$ Then $X$ is a positive multiple of $\ell^2$ (in particular, this implies that $X$ is Euclidean). 
\end{lemma}
\begin{proof}

For brevity, we will omit some details. For any $u\in \mathbb R^2$ and $t\in [0,1]$, let $D_t(u)$ be the disk of diameter $[-tu,u].$ If $u\neq 0$, the reader should confirm that $v\in D_0(u)^\circ$ if and only if $\ell_v$ intersects $(0,u)$. (This can be proven via the inscribed angle theorem, or
algebraic manipulation of the inner product). Suppose for the sake of contradiction that we have $u\in B_X$ and $v\in D_0(u)^\circ\setminus B_X^\circ.$ Then there exists some $w\in \ell_v\cap (0,u)$. Since $||v||_X\geq 1$, we can define $$\widehat v=\frac{v}{\,||v||_X\!\!
\!}\in \partial B_X\qquad\text{and}\qquad \widehat w=\frac{w}{\,||v||_X\!\!\!}\in (0,u).$$ Then $\widehat w\in \ell_{\widehat v}\cap (0,u)\subset \ell_{\widehat v}\cap B_X^\circ$, which contradicts the assumption in the statement of the lemma. Therefore, we have $u\in B_X\Longrightarrow D_0(u)^\circ\subset B_X^\circ.$ Since $B_X$ is closed, this gives $u\in B_X\Longrightarrow D_0(u)\subset B_X$. From here, we get a well-defined supremum
\begin{equation}\label{circle-sup-eq}
m=\sup\{t\in [0,1]:u\in B_X\Longrightarrow D_t(u)\subset B_X\}.    
\end{equation}
Then $m\in [0,1]$ and since $B_X$ is closed, we can see that $u\in B_X\Longrightarrow D_m(u)\subset B_X.$ Since $B_X$ is compact, there exists $u\in B_X$ with maximal $\ell_2$-norm. If $m=1,$ then $$B_X\subset |u| \cdot B_{\ell^2}= D_1(u)\subset B_X.$$ Thus $B_X=|u|\cdot B_{\ell^2}$ and hence $X= |u|\cdot \ell^2$ is a positive multiple of $\ell^2$ (if $m=1$).

\begin{figure}
    \centering
        \vspace*{5mm}
    \begin{tikzpicture}[scale=3]
    
    \begin{scope}[shift={(-0.3,-0.2)}]
    \draw[domain=0:360,samples=500,thick] plot (\x:{(11/9+sin(\x))/(11/9+1)});
    \foreach \i in {0,12,...,348}{
    \pgfmathsetmacro{\x}{(0.9*sin(\i)-sqrt(0.81*sin(\i)*sin(\i)+0.4)}
    \draw[scale=0.25] ({\x*cos(\i)},{\x*sin(\i)}) circle (\x);
    }
    
    \draw (0.74,0.3) node{$L$};
    
    \end{scope}
    
    \draw (-0.3,-0.6) node{(a)};
    \draw (2,-0.6) node{(b)};

    \begin{scope}[shift={(2,0)}]
    \draw[domain=0:360,samples=500,thick] plot (\x:{(11/9+sin(\x))/(11/9+1)});
     \draw[domain=-30:210,samples=500] plot (\x:{(11/9-sin(\x))/(11/9+1)});
     \draw[domain=-30:-42,samples=500,dash pattern=on 2 off 2] plot (\x:{(11/9-sin(\x))/(11/9+1)});
     \draw[domain=210:222,samples=500,dash pattern=on 2 off 2] plot (\x:{(11/9-sin(\x))/(11/9+1)});
    \draw (0,0.45) circle(0.55);
    \draw[] (0,0.5-0.3025/2) circle (0.5+0.3025/2);

    \draw[thick,->] (0,0) node{$\bullet$} -- (0,1);
    
    \draw (-0.07,0.5) node{$u$};
    
    \draw (0.34,0.48) node{$D_m(u)$};
    \draw (0.32,-0.36) node{$D_s(u)$};
    \draw (0.74,0.3) node{$L$};
    
    \draw (-0.8,-0.32) node{$-L$};
    
    \end{scope}
    
    \end{tikzpicture}
    \caption{Sweeping out a lima\c{c}on to get a larger disk}
    \label{fig:limacon}
\end{figure}
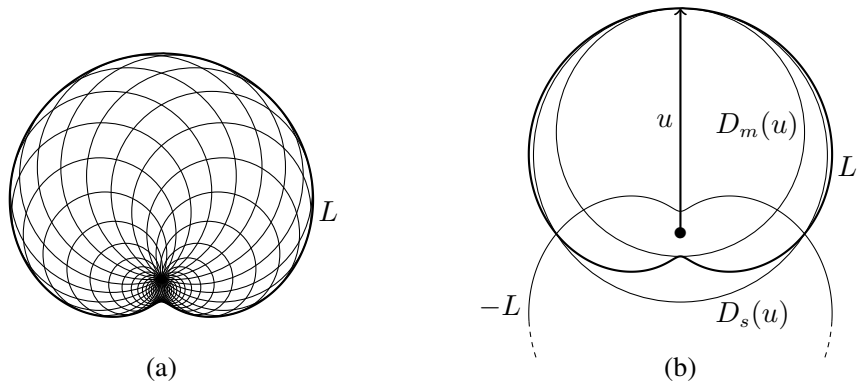

To complete the proof, we will assume that $m<1$ and derive a contradiction. Let $$a=\frac{1+m}{1-m}\qquad\text{and}\qquad s=\left(\frac{a}{a+1}\right)^{\!\!2}=\left(\frac{m+1}{2}\right)^{\!\!2}.$$ The reader may check that $m<s<1.$ We will prove $u\in B_X\Longrightarrow D_{s}(u)\subset B_X$, which contradicts the definition of $m$ in (\ref{circle-sup-eq}). Towards this end, suppose that $u\in B_X$ and $u\neq 0$ (we clearly have  $D_s(0)=\{0\}\subset B_X$). Let $r_0=|u|$ and $\theta_0=\arg(u)$ (equivalently, we could write $u=r_0\exp(i\theta_0)$ in complex notation). We define sets
\begin{equation}\label{limacon-def}
    {\pm}L=\left\{(r,\theta)\in \mathbb R^2:\frac{r}{r_0}\leq \frac{a\pm \cos (\theta-\theta_0)}{a+1} \right\},
\end{equation} where $(r,\theta)$ denotes polar coordinates. The set $L$ is a filled lima\c{c}on, which is linearly equivalent to the more familiar form $\{(r,\theta)\in \mathbb R^2:r\leq a+\cos\theta\}.$ Then we have
$$L=\!\!\!\bigcup_{v\in \partial D_m(u)}\!\!\!D_0(v).$$ (This fact, illustrated in Figure \ref{fig:limacon}(a), is commonly expressed by saying that the lima\c{c}on $\partial L$ is the ``envelope'' of $\{\partial D_0(v):v\in \partial D_m(u)\}$ \cite{lockwood}.) Since $v\in \partial D_m(u)\subset B_X$ implies $D_0(v)\subset B_X$, we see that $L\subset B_X$ and thus $-L\subset -B_X=B_X.$ We have
\begin{equation}\label{polar-circle}
  D_s(u)=\left\{(r,\theta)\in \mathbb R^2:\frac{r}{r_0}\leq f\big(\cos(\theta-\theta_0)\big)\right\}  
\end{equation}
 in polar coordinates $(r,\theta),$ where $2f(t)=(1-s)t+\sqrt{4s+(1-s)^2t^2}.$ Then
\begin{equation}\label{limacon-containment}
t\in [0,1]\Longrightarrow f(t)\leq \frac{a+t}{a+1}\quad\text{and}\quad t\in [-1,0]\Longrightarrow f(t)\leq \frac{a-t}{a+1}.   
\end{equation}
Considering the sign of $\cos(\theta-\theta_0)$, we may compare (\ref{limacon-def}) and (\ref{polar-circle}) via (\ref{limacon-containment}), to see that: 
\begin{align*}
    &\text{if }\theta -\theta_0\in \left[\tfrac{-\pi}{2},\tfrac{\pi}{2}\right],\text{ then } (r,\theta)\in D_s(u)\Longrightarrow (r,\theta)\in L\subset B_X,\\
    &\text{if }\theta -\theta_0\in \left[\tfrac{\pi}{2}\,,\tfrac{3\pi}{2}\right],\text{ then } (r,\theta)\in D_s(u)\Longrightarrow (r,\theta)\in -L\subset B_X.
\end{align*}
Therefore $D_s(u)\subset B_X$, as desired. This argument is illustrated in Figure \ref{fig:limacon}(b).
\end{proof}

\subsection{Quarter-turn symmetry} We now return to $\varpi$, to recount and expand upon the result of \cite{pi-values} on norms with quarter-turn symmetry, which states that $\varpi(X)\geq \pi$ whenever $iX=X$ (this means that pushing forward by a quarter-turn does not change lengths). Recall that $iX=X\Longleftrightarrow iB_X=B_X$ (i.e.~$B_X$ has quarter-turn symmetry). We begin with a crucial lemma, which compares lengths under such a norm to Euclidean angles.

\begin{lemma}\label{quarter-lemma} Fix any norm $X\in \mathcal M$ with $iX=X$.
\begin{enumerate}
    \item[(a)] Suppose that $v\neq 0$ and $p,q\in \ell_v$ are distinct. If $\theta$ denotes the angle measure between the vectors $p$ and $q$, then we have $d_X(p,q)>\theta\cdot ||v||_X.$ 
    \item[(b)] If $\varphi$ is a path along $\partial B_X$ that sweeps out an angle of $\theta$ (centered at the origin), then $\textup{len}_X\varphi\geq \theta.$ If $\varphi$ is a polygonal path (and not just a point), then $\textup{len}_X\varphi>\theta.$
\end{enumerate}
\end{lemma}
\begin{proof}
The underlying argument in this proof follows \cite{pi-values} extremely closely.

\smallskip (a) We first consider the Euclidean geometry. Define $r=|p|,$ $s=|q|,$ $b=|p-q|$ and $h=|v|,$ as shown in Figure \ref{fig:quarter-turn-diagram}(a). We have $bh=rs\sin \theta$, because both sides equal twice the area of $\triangle 0pq.$ By the law of cosines in the same triangle, we also have $b^2=r^2+s^2-2rs\cos \theta.$ Note that $t+1/t\geq 2$ for all $t>0$. Therefore, we have
$$\frac{b}{h}=\frac{b^2}{bh}=\frac{r^2+s^2-2rs\cos\theta}{rs\sin \theta}=\frac{r/s+s/r-2\cos \theta}{\sin\theta} \geq \frac{2-2\cos \theta}{\sin\theta}.$$ Since $p$ and $q$ are not parallel, we have $\theta<\pi$. The tangent half-angle formula gives
$$\frac{b}{h}\geq \frac{2-2\cos \theta}{\sin \theta}=2\tan (\theta /2)>\theta,$$ where the last inequality follows because $\tan t>t$ for all $0<t<\pi/2.$ 

Now, we must consider length in terms of $X.$ Since $p,q\in \ell_v$, we have $p-q=tiv$ for some $t\in \mathbb R.$ For any norm $Y\in \mathcal  M$ with $iY=Y$ (e.g.~$X$ or $\ell^2$), it follows that $$||p-q||_Y=||tiv||_Y=|t|\cdot ||iv||_Y=|t|\cdot ||v||_Y.$$ This allows us to translate the result for $\ell^2$ to any norm $X$ with $iX=X$:
$$||p-q||_X=|t|\cdot ||v||_X=\frac{|p-q|}{|v|}\cdot||v||_X=\frac{b}{h}\cdot||v||_X>\theta\cdot ||v||_X.$$

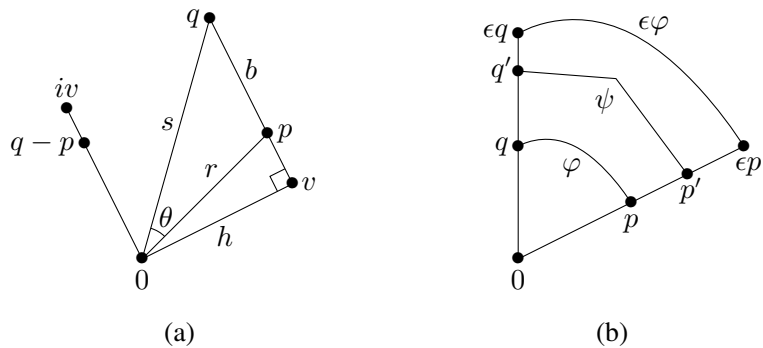
\begin{figure}
    \centering
        \vspace*{5mm}
    \begin{tikzpicture}
    
    \begin{scope}
    \draw (-1,2) node{$\bullet$} -- (0,0) node{$\bullet$} -- (2,1) node{$\bullet$};
    
    \draw ({0.9-5/3},{3.2-5/3}) node{$\bullet$};
    
    \draw ({5/3},{5/3}) node{$\bullet$} -- (0,0) -- (0.9,3.2) node{$\bullet$} -- cycle;
    
    \draw ({5/3},{5/3}) -- (2,1);
    
    \draw (0.9*2,0.9*1) -- (1.8-0.1,0.9+0.2) -- (2-0.1,1+0.2);
    
    \draw (0.3,0.3) arc (45:74:{0.3*sqrt(2)});
    
    \draw (0.32,0.56) node{$\theta$};
    
    \draw (2.22,1) node{$v$};
    \draw (-1,2.28) node{$iv$};
    \draw (-1.33,1.53) node{$q-p$};
    \draw (0,-0.27) node{$0$};
    \draw (0.68,3.2) node{$q$};
    \draw (1.9,1.67) node{$p$};
    
    \draw (1.1,0.3) node{$h$};
    \draw (1.45,2.53) node{$b$};
    \draw (0.34,1.8) node{$s$};
    \draw (0.92,1.18) node{$r$};
    
    \draw (0.5,-1) node{(a)};
    \end{scope}
    
    \begin{scope}[shift={(5,0)}]
    
    \draw (1.25,-1) node{(b)};
    
    \draw plot [smooth, tension=1] coordinates { (1.5,0.75) (0.75,1.5) (0,1.5)};
        \draw plot [smooth, tension=1] coordinates { (3,1.5) (1.5,3) (0,3)};
        
        \draw (0,3) node{$\bullet$} -- (0,1.5) node{$\bullet$} -- (0,0) node{$\bullet$} -- (1.5,0.75) node{$\bullet$} -- (3,1.5) node{$\bullet$};
        \draw (0,2.5) node{$\bullet$} -- (1.3,2.4) -- (2.25,1.125) node{$\bullet$};
        
        \draw (0,-0.27) node{$0$};
        \draw (0.71,1.21) node{$\varphi$};
        \draw (1.8,3.1) node{$\epsilon\varphi$};
        \draw (1.15,2.1) node{$\psi$};
    
        \draw (-0.22,1.5) node{$q$};
        \draw (-0.22,2.52) node{$q'$};
        \draw (-0.3,3) node{$\epsilon q$};
    
        \draw (1.5,0.47) node{$p$};
        \draw (2.3,0.91) node{$p'$};
        \draw (3.06,1.23) node{$\epsilon p$};
    
    \end{scope}
    
    \end{tikzpicture}
    
    \caption{Lengths in terms of a norm $X$ such that $iX=X$}
    \label{fig:quarter-turn-diagram}
\end{figure}

\smallskip (b) First suppose that $\varphi=[p,q]$ with $p\neq q.$ Let $\ell$ denote the line through $p$ and $q$. Since $[p,q]\subset \partial B_X$, Lemma \ref{extreme-pts}(a) states that $B_X\cap \ell$ is a face of $B_X.$ Because $B_X\cap \ell$ is a proper face of $B_X$, we have $B_X\cap \ell\subset \partial B_X$ and thus $\ell\cap B_X^\circ=\emptyset.$ In particular, this shows that $0\notin \ell,$ so we can write $\ell=\ell_v$ for some $v\in \ell$ ($v$ is the unique point along $\ell$ of minimal $\ell^2$-norm). Since $v\notin B^\circ_X$, we have $||v||_X\geq 1.$ Therefore (1) gives
$$\text{len}_X\varphi=d_X(p,q)>\theta\cdot ||v||_X\geq \theta.$$ This proves the desired result whenever $\varphi$ is a line segment. Since lengths are additive under concatenation of paths, we also get the desired result for any polygonal path. 

In proving the general case, we may assume that $\theta<\pi$ (by the additivity of length). Then if $\varphi$ goes from $p$ to $q$, the angle $\theta$ swept out by $\varphi$ is the angle between the vectors $p$ and $q$ (which is unchanged if $p$ or $q$ is scaled by a positive number). Fix some $\epsilon>1$. Let $P\subset \mathbb R^2$ be a convex polygon with $iP=P$ and $B_X\subset P\subset \epsilon B_X$. (To prove that $P$ exists, we proceed as follows. Let $x_0,x_1,x_2,\dots \subset\partial(\epsilon B_X)$ be a dense sequence. Then $B_X\subset\text{conv}\{x_n:n\in \mathbb N\}^\circ$ and thus the sets $\big\{\text{conv}\{x_0,x_1,\dots,x_n\}^\circ:n\in \mathbb N\big\}$ form an open cover of $B_X.$ By compactness, one of these sets actually contains $B_X.$ We then define $P=\text{conv}\{ax_k:k=0,1,\dots,n\text{ and }a=\pm 1 \text{ or }\pm i\}$.) Then there exist unique points $p',q'\in \partial P$ such that $p'=sp$ and $q'=tq$ for some $s,t\in [1,\epsilon].$ Let $\psi$ denote the (shorter) path along $\partial P$ from $p'$ to $q'.$ Notice that $\psi$ also sweeps out an angle of $\theta$ and $P$ is a polygon, so the above case gives $\text{len}_Y\psi>\theta$, where $Y\in \mathcal M$ is the unique norm with $B_Y=P.$ Then $\psi\prec [p',\epsilon p]\bullet \epsilon\varphi\bullet [\epsilon q,q']$, where the latter is a convex path because it is a part of the boundary of $\{tv:t\in [0,\epsilon]\text{ and }v\in \varphi\}.$ These convex paths are illustrated in Figure \ref{fig:quarter-turn-diagram}(b). Since $p'\in [p,\epsilon p],$ we have $$d_Y(p',\epsilon p)\leq d_Y(p,\epsilon p)=(\epsilon-1)||p||_Y.$$ Analogous reasoning gives $d_Y(q',\epsilon q)\leq (\epsilon-1)||q||_Y.$ By Lemma \ref{convex-paths}(a), we then have
\begin{align*}
    \theta<\text{len}_Y\psi &\leq \text{len}_Y\big([p',\epsilon p]\bullet \epsilon\varphi\bullet [\epsilon q,q']\big)\\
    &=d_Y(p',\epsilon p)+d_Y(q',\epsilon q)+\text{len}_Y(\epsilon \varphi)\\
    &\leq (\epsilon-1)\big(||p||_Y+||q||_Y\big)+ \text{len}_Y(\epsilon\varphi).
\end{align*}
Since $B_X\subset P=B_Y$, we have $\text{len}_X\geq \text{len}_Y$ and $||\cdot ||_X\geq ||\cdot||_Y.$ Therefore,
\begin{align*}
    \epsilon \cdot\text{len}_X\varphi&\geq \epsilon\cdot \text{len}_Y\varphi=\text{len}_Y(\epsilon \varphi)\\
    &>\theta -(\epsilon-1)\big(||p||_Y+||q||_Y\big)\\
    &\geq \theta -(\epsilon -1)\big(||p||_X+||q||_X\big).
\end{align*}
Taking the limit as $\epsilon\rightarrow 1$ yields the desired inequality $\text{len}_X\varphi\geq \theta$.\end{proof}

This lemma carries most of the burden of proving that $\varpi(X)\geq \pi$ when $iX=X$, as well as classification of the equality case. But the condition $iX=X$ is not ideal, in that it is not preserved under linear equivalence (consider ellipses).  Thus we define 
$$\mathcal Q=\{Y\in \mathcal M:Y\text{ is linearly equivalent to some }X\in \mathcal M\text{ with }iX=X\}.$$ Then $\mathcal Q$ is obviously closed under linear equivalence, so it provides a ``coordinate-free'' notion of norms with quarter-turn symmetry.\footnote{A more coordinate-free definition is $X\in \mathcal Q$ if and only if $SX=X$ for some $S\in \text{GL}(2,\mathbb R)$ of order 4.} With this notion in hand, the quarter-turn symmetry results can be extended to $\mathcal Q$, giving a characterization of $\mathcal E$ in terms of $\varpi$.

\begin{proposition}\label{quarter-turn-prop} 
\begin{enumerate}
    \item[(a)] For any norm $X\in \mathcal M$ such that $iX=X$, we have $\varpi(X)\geq \pi.$ Moreover, we have $\varpi(X)=\pi$ if and only if $X$ is a positive multiple of $\ell^2$.
    \item[(b)] We have $\varpi(\mathcal Q)=[\pi,4]$ and $\mathcal E=\{X\in \mathcal Q:\varpi(X)=\pi\}$.
\end{enumerate}
\end{proposition}

\begin{proof}
(a) Note that $\varpi(X)\geq \pi$ follows at once from Lemma \ref{quarter-lemma}(b), taking $\varphi=\partial B_X$, which sweeps out a full angle of $\theta=2\pi$. If $X$ is a positive multiple of $\ell^2$, then $X\in \mathcal E$ and thus $\varpi(X)=\pi.$ Conversely, now suppose that $X$ is not a positive multiple of $\ell^2.$ The contrapositive of Lemma \ref{classify-circles} states that there is some $v\in \partial B_X$ with $B_X^\circ \cap \ell_v\neq \emptyset.$ Because the line $\ell_v$ intersects the interior of $B_X$, it must intersect the boundary $\partial B_X$ in at least two points, so there is some $u\in \ell_v\cap \partial B_X$ with $u\neq v.$ Let $\varphi$ (resp.~$\psi$), denote the shorter (resp.~longer) path along $\partial B_X$ from $u$ to $v.$ If $\theta$ is the angle formed by $u$ and $v$, then $\varphi$ (resp.~$\psi$) sweeps out an angle of $\theta$ (resp.~$2\pi-\theta$). It follows that
$$\text{len}_X(\partial B_X)=\text{len}_X\psi+\text{len}_X\varphi\geq 2\pi-\theta+d_X(p,q)>2\pi-\theta+\theta=2\pi,$$ by Lemmas \ref{quarter-lemma}(a) and (b). Therefore, we get the \textit{strict} inequality $\varpi(X)>\pi$.

\smallskip (b) If $Y\in \mathcal Q$, then $Y$ is linearly equivalent to some $X\in \mathcal M$ such that $iX=X.$ Then $\varpi(Y)=\varpi(X)\geq \pi$ by Lemma \ref{linear-equiv} and (a). Since the upper bound of $\varpi(Y)\leq 4$ was already established in Proposition \ref{upper-bound}, we have $\varpi(\mathcal Q)\subset [\pi,4].$ To show the reverse, note that $\ell^p$ has quarter-turn symmetry for all $p\in [1,\infty].$ We know that $\varpi(\ell^2)=\pi$ and $\varpi(\ell^1)=4$, so the continuity of $p\mapsto \ell^p\mapsto \varpi(\ell^p)$ (with the right metric on $\mathcal M$, which we will not discuss in detail here) and the intermediate value theorem give $$[\pi,4]\subset\varpi\big(\{\ell^p:1\leq p\leq 2\}\big)\subset \varpi(\mathcal Q).$$ For more details on the function $p\mapsto \varpi(\ell^p),$ including a proof of continuity, see \cite{lp-norms}.

We already know that $\ell^2\in \mathcal Q$ and $\varpi(\ell^2)=\pi$. But since $\{X\in \mathcal Q:\varpi(X)=\pi\}$ is closed under linear equivalence, this implies that $\mathcal E\subset \{X\in \mathcal Q:\varpi(X)=\pi\}.$ Conversely, consider $Y\in \mathcal Q$ with $\varpi(Y)=\pi.$ Then $Y$ is linearly equivalent to some $Z\in \mathcal M$ with $iZ=Z.$ We have $\varpi(Z)=\varpi(Y)=\pi$ by Lemma \ref{linear-equiv}, so $Z$ is a positive multiple of $\ell^2$ by (a). Hence, $Z$ is Euclidean and therefore $Y$ is Euclidean as well. \end{proof}

\begin{acknowledgment}{Acknowledgement.}
I would like to extend my thanks to Cornelia Van Cott for the talk that inspired this article and the encouragement that helped it come to fruition. 
\end{acknowledgment}

\begin{biog}
\item[Nikhil Henry Bukowski Sahoo] hails from the East Bay, where he attended the Peralta Community
Colleges and where he eventually earned a B.A. in mathematics from the University of California, Berkeley.
He is currently a Ph.D. student in mathematics at Cornell University studying equivariant symplectic topology.
\begin{affil}
Cornell University, Ithaca NY 14850\\
nhs58@cornell.edu
\end{affil}

\end{biog}

\vfill\eject

\end{document}